\documentclass[10pt]{amsart}
\usepackage{latexsym,amscd,amssymb,verbatim} 

\theoremstyle{plain}
  \newtheorem{theorem}{Theorem}[]
  \newtheorem{proposition}[theorem]{Proposition}
  \newtheorem{lemma}[theorem]{Lemma}

\theoremstyle{definition}

  \newtheorem{question}[theorem]{Question}
 \theoremstyle{remark}

\numberwithin{equation}{section}

\def\ZZ{{\mathbb Z}}

\def\Spec{{\operatorname{Spec}}}

\def\aa{{\mathbf{a}}}
\def\bb{{\mathbf{b}}}
\def\AA{{\mathbf{A}}}
\def\BB{{\mathbf{B}}}

\def\XX{{X}}

\def\gammma{{c}}

\def\SSS{{\mathcal{S}}}

\def\mm{{\mathfrak{m}}}

\begin{document}

\title[A universal coefficient theorem for Gau{\ss}'s Lemma]
{A universal coefficient theorem for Gau{\ss}'s Lemma}

\author{William Messing}
\email{messing@math.umn.edu}
\author{Victor Reiner}
\email{reiner@math.umn.edu}

\address{School of Mathematics\\
University of Minnesota\\
Minneapolis, MN 55455}

\dedicatory{To J\"urgen Herzog on his 70th birthday}

\thanks{Second author partially supported by NSF grant DMS-1001933.}


%




\begin{abstract}
We shall prove a version of Gau\ss's Lemma. 
It works in $\ZZ[\aa,\AA,\bb,\BB]$ where 
$\aa=\{a_i\}_{i=0}^m,
\AA=\{A_i\}_{i=0}^m,
\bb=\{b_i\}_{j=0}^n,
\BB=\{B_j\}_{j=0}^n,
$
and constructs polynomials $\{ \gammma_k \}_{k=0,\ldots,m+n}$
of degree at most $\binom{m+n}{n}$ in each variable set
$\aa,\AA,\bb,\BB$, with this property: setting
$$
\sum_k C_k \XX^k =
\sum_{i} A_i \XX^i
\cdot
\sum_{j} B_j \XX^j ,
$$
for elements $a_i,A_i,b_j,B_j$ in any
commutative ring $R$ satisfying
$$
1=\sum_i a_i A_i =\sum_j b_j B_j,
$$
the elements $c_k=c_k(a_i,A_i,b_j,B_j)$ will satisfy
$1=\sum_k \gammma_k C_k$.
\end{abstract}

\maketitle

\section{The statement}
\label{statement-section}

Let $R$ be a commutative ring.  Consider two elements 
$A(\XX)=\sum_{i=0}^m A_i \XX^i$ and $B(\XX)=\sum_{j=0}^n B_j \XX^j$ 
in $R[\XX]$, with product 
$C(\XX)=A(\XX)B(\XX) = \sum_{k=0}^{m+n} C_k \XX^k$,
so that 
$
C_k=\sum_{i+j=k} A_i B_j.
$
A version of Gau\ss's Lemma, called {\it Gau\ss-Joyal de pauvre} 
in \cite[\S II Lemme 2.6]{LombardiQuitte}, asserts the following.

\begin{proposition} 
\label{Gauss's-lemma}
If both $A(\XX), B(\XX)$ have the property that their coefficient sequences
generate the unit ideal $R$, then the same is true of their product $C(\XX)$,
that is,
$
(A_0,\ldots,A_m)=R=(B_0,\ldots,B_n)
$
implies $(C_0,\ldots,C_{m+n})=R$.
\end{proposition}

\noindent
Standard proofs of this appeal to Zorn's Lemma or some weaker version,
to show that if the ideal $(C_0,\ldots,C_{m+n})$ is not $R$,
then it is contained in some maximal ideal $\mm$ of $R$.
This leads to the contradiction in the integral domain $R/\mm[\XX]$,
that $A(\XX), B(\XX)$ represent nonzero elements
but their product represents zero.  Gau\ss's Lemma has several proofs
avoiding any version of Zorn's Lemma;  we discuss some of these 
in Section~\ref{historical-section} below.

Our goal here, however, is a construction
of ``universal'' coefficients $\gammma_0,\ldots,\gammma_{m+n}$ satisfying 
$1=\sum_{k=0}^{m+n} \gammma_k C_k$, expressed in terms of any
$a_i, b_j$ in $R$ that satisfy
\begin{equation}
\label{gcd-equation}
1=\sum_{i=0}^{m} a_i A_i=\sum_{j=0}^n b_j B_j.
\end{equation}
To do this, we work in a polynomial algebra, 
$$
S=\ZZ[\aa,\bb,\AA,\BB]
=\ZZ[a_0,\ldots,a_m,b_0,\ldots,b_n,A_0,\ldots,A_m,B_0,\ldots,B_n]
$$
and consider the elements 
\begin{equation}
\label{generators-of-relation-ideal}
a:=1-\sum_{i=0}^{m} a_i A_i,\qquad b:=1-\sum_{j=0}^n b_j B_j, \qquad 
C_k:=\sum_{i+j=k} A_i B_j, 
\end{equation}
for $k=0,1,\ldots,m+n$.  We shall prove the following equivalent of Gau\ss's Lemma.
\begin{proposition}
\label{constructive-Gauss's-lemma}
There exist in $S$ polynomials 
$\alpha, \beta, \gammma_0,\gammma_1,\ldots,\gammma_{m+n}$ expressing
\begin{equation}
\label{unit-ideal-expression}
1 = \alpha a + \beta b + \sum_{k=0}^{m+n} \gammma_k C_k.
\end{equation}
Consequently, for any commutative ring $R$ and elements satisfying
\eqref{gcd-equation}, the elements $\gammma_k=\gammma_k(\aa,\bb,\AA,\BB)$ will satisfy
$
1= \sum_{k=0}^{m+n} \gammma_k C_k
$
in $R$.
\end{proposition}

\section{A proof which is not explicit}
\label{less-constructive-section}

We first give a proof of the existence of polynomials 
$\alpha,\beta,\gammma_k$ as in Proposition~\ref{constructive-Gauss's-lemma},
via induction on $m+n$. 
In Section~\ref{more-constructive-section} below, we will reinterpret
it to give explicit recursive formulas for $\alpha,\beta,\gammma_k$.
These formulas also lead to a simple bound (Proposition~\ref{degree-bound-proposition} below) 
on the degrees of 
$\alpha,\beta,\gammma_k$ when considered as polynomials in 
each of the variable sets $\aa,\bb,\AA,\BB$.

For the base case of the induction, let $m=n=0$.
One then checks directly that 
$\alpha=1,\beta=a_0A_0$, and $\gammma_0=a_0b_0$ suffice:
$$
1= 1 \cdot (1-a_0 A_0) + a_0 A_0 \cdot (1- b_0 B_0) + a_0 b_0 \cdot A_0 B_0.
$$

In the induction on $m+n$, it will be important 
to emphasize the dependence of various objects on $m,n$.
Changing notation, denote by $S^{(m,n)}$ the polynomial ring $S$, 
and denote by $a^{(m)},b^{(n)},C^{(m,n)}_k$ the elements $a,b,C_k$ 
appearing in \eqref{generators-of-relation-ideal}.  
Let $Q^{(m,n)}$ denote the quotient ring of $S^{(m,n)}$ by the ideal generated
by these elements.  
Proposition~\ref{constructive-Gauss's-lemma} 
then asserts $1=0$ in $Q^{(m,n)}$, that is $Q^{(m,n)}$ is the 
zero ring.

One easily checks the following comparisons for $m,n \geq 1$
\begin{equation}
\label{comparisons}
\begin{aligned}
a^{(m-1)} &= a^{(m)}+a_m A_m,\\
b^{(n-1)} &= b^{(n)}+b_n B_n, \\
C^{(m-1,n)}_k &= 
 \begin{cases} 
   C^{(m,n)}_k & \text{ if } k < m,\\
   C^{(m,n)}_k - A_m B_{k-m} & \text{ if } k \geq m, \\
 \end{cases} \\
C^{(m,n-1)}_k &= 
 \begin{cases} 
   C^{(m,n)}_k & \text{ if } k < n,\\
   C^{(m,n)}_k - A_{k-n} B_{n} & \text{ if } k \geq n. \\
 \end{cases}
\end{aligned}
\end{equation}
\noindent
which show that the principal ideals $(A_m),(B_n)$ satisfy\footnote{These isomorphisms hold even when $m=0$ or $n=0$,
provided one adopts the convention that
$Q^{(-1,n)}=Q^{(m,-1)}=0$:  the vanishing of $a^{(0)}=1-a_0 A_0$
in $Q^{(0,n)}$ makes $A_0$ a unit, and hence $Q^{(0,n)}/(A_0)=0$.}
\begin{equation}
\label{quotient-recursions}
\begin{aligned}
Q^{(m,n)}/(A_m) &\cong Q^{(m-1,n)} \text{ for }m \geq 1,\\
Q^{(m,n)}/(B_n) &\cong Q^{(m,n-1)} \text{ for }n \geq 1.\\
\end{aligned}
\end{equation}

Now assuming that $m+n \geq 1$, one 
has $Q^{(m,n)}/(B_n) \cong Q^{(m,n-1)} =0$ using
\eqref{quotient-recursions} and induction.
Hence it suffices to show that the ideal $(B_n)=0$.
Note that multiplication by $B_n$ gives a
surjection  
$
Q^{(m,n)} \twoheadrightarrow (B_n)
$
that factors through $Q^{(m,n)}/(A_m)$, since
$A_m$ annihilates $B_n$.  
As $Q^{(m,n)}/(A_m) \cong Q^{(m-1,n)} =0$,
by \eqref{quotient-recursions} and induction,
this shows $(B_n)=0$, completing the proof.

\section{The explicit recursions}
\label{more-constructive-section}

Here we recursively produce,
for $m,n \geq 0$, polynomials $\alpha^{(m,n)},\beta^{(m,n)},\gammma^{(m,n)}_k$
in the ring $S=S^{(m,n)}$ with the property that
\begin{equation}
\label{desired-constructive-equality}
1=\alpha^{(m,n)} a^{(m)} + \beta^{(m,n)} b^{(n)}
   + \sum_{k=0}^{m+n} \gammma^{(m,n)}_k C^{(m,n)}_k.
\end{equation}
\noindent
One easily checks that when $m=0$ one can choose:
\begin{equation}
\label{initial-conditions-for-zero-m}
\alpha^{(0,n)}=1, \,\,
\beta^{(0,n)}=a_0 A_0, \,\,
\gammma^{(0,n)}_k=a_0 b_k.
\end{equation}
By symmetry, when $n=0$ one can choose:
\begin{equation}
\label{initial-conditions-for-zero-n}
\alpha^{(m,0)}=b_0 B_0, \,\,
\beta^{(m,0)}=1, \,\,
\gammma^{(m,0)}_k=a_k b_0.
\end{equation}
For definiteness, when $m=n=0$,
we choose to use \eqref{initial-conditions-for-zero-m}
instead of \eqref{initial-conditions-for-zero-n}.

Now assume that $m, n>0$.  By induction on $m+n$, assume that one 
has constructed 
$\alpha^{(m,n-1)},\beta^{(m,n-1)},\gammma_k^{(m,n-1)},
\alpha^{(m-1,n)},\beta^{(m-1,n)},\gammma_k^{(m,n-1)}$ 
satisfying:
$$
\begin{aligned}
1&=\alpha^{(m,n-1)} a^{(m)} + \beta^{(m,n-1)} b^{(n-1)} 
   + \sum_{k=0}^{m+n-1} \gammma^{(m,n-1)}_k C^{(m,n-1)}_k, \\
1&=\alpha^{(m-1,n)} a^{(m-1)} + \beta^{(m-1,n)} b^{(n)} 
   + \sum_{k=0}^{m+n-1} \gammma^{(m-1,n)}_k C^{(m-1,n)}_k.
\end{aligned}
$$
Using  \eqref{comparisons}, these become
\begin{equation}
\label{first-almost-the-identity}
1 = \alpha^{(m,n-1)} a^{(m)} + \beta^{(m,n-1)} b^{(n)}
 + \sum_{k=0}^{m+n-1} \gammma^{(m,n-1)}_k C^{(m,n)}_k 
 + B_n d 
\end{equation}
\begin{equation}
\label{second-almost-the-identity}
1 = \alpha^{(m-1,n)} a^{(m)} + \beta^{(m-1,n)} b^{(n)}
 + \sum_{k=0}^{m+n-1} \gammma^{(m-1,n)}_k C^{(m,n)}_k 
 + A_m e.
\end{equation}
where we have used two auxiliary polynomials
\begin{equation}
\label{auxiliary-polynomials}
\begin{aligned}
d&=
b_n \beta^{(m,n-1)} - \sum_{k=n}^{m+n-1} \gammma^{(m,n-1)}_k A_{k-n}, \\
e&=
a_m \alpha^{(m-1,n)} - \sum_{k=m}^{m+n-1} \gammma^{(m-1,n)}_k B_{k-m}. \\
\end{aligned}
\end{equation}
\noindent
Now one can use \eqref{second-almost-the-identity} 
to replace $d$ in \eqref{first-almost-the-identity} by
$$
d=d \cdot 1 = d \left( 
\alpha^{(m-1,n)} a^{(m)} + \beta^{(m-1,n)} b^{(n)}
 + \sum_{k=0}^{m+n-1} \gammma^{(m-1,n)}_k C^{(m,n)}_k 
 + A_m e
\right).
$$
This yields the following expression:
$$
\begin{aligned}
1 &=\left( \alpha^{(m,n-1)}+B_n d \alpha^{(m-1,n)} \right) a^{(m)}
 + \left( \beta^{(m,n-1)}+B_n d \beta^{(m-1,n)} \right) b^{(n)} \\
&+ \sum_{k=0}^{m+n-1} \left( \gammma^{(m,n-1)}_k+B_n d \gammma^{(m-1,n)}_k \right) C^{(m,n)}_k
 + d e \cdot A_m B_n.
\end{aligned}
$$
Thus, if one recursively defines
\begin{equation}
\label{alpha-beta-gamma-recursion}
\begin{aligned}
\alpha^{(m,n)}&=\alpha^{(m,n-1)}+B_n d \alpha^{(m-1,n)},\\
\beta^{(m,n)}&=\beta^{(m,n-1)}+B_n d \beta^{(m-1,n)},\\
\gammma^{(m,n)}_k&=\gammma^{(m,n-1)}_k+B_n d \gammma^{(m-1,n)}_k
\text{ for } k=0,1,\ldots,m+n-1\\
\gammma^{(m,n)}_{m+n}&=d e \quad
   \text{ (since } C^{(m,n)}_{m+n}=A_m B_n\text{)},
\end{aligned}
\end{equation}
then one obtains coefficients 
satisfying \eqref{desired-constructive-equality}.
We have thus proven the following.

\begin{proposition}
\label{summary-prop}
The polynomials $\alpha^{(m,n)}, \beta^{(m,n)}, \gammma_k^{(m,n)}$
defined by recursions 
\eqref{alpha-beta-gamma-recursion} and \eqref{auxiliary-polynomials},
with initial conditions given by 
\eqref{initial-conditions-for-zero-m} and 
\eqref{initial-conditions-for-zero-n},
satisfy \eqref{desired-constructive-equality} 
(and \eqref{unit-ideal-expression}). 
\end{proposition}

\begin{proposition}
\label{degree-bound-proposition}
The polynomials $\alpha^{(m,n)}, \beta^{(m,n)}, \gammma_k^{(m,n)}$
in Proposition~\ref{summary-prop} have degree bounded 
by $\binom{m+n}{m}$ in each variable
set $\aa,\bb,\AA,\BB$.
\end{proposition}
To establish this, we use induction on $m+n$ to 
prove a slightly more precise statement.  Let
$
N:=\binom{m+n}{m},
N':=\binom{m+n-1}{m},
N'':=\binom{m+n-1}{m-1}
$
so that $N=N'+N''$.  It is straightforward to check that
the bounds on the total degree in the following
table are valid in each variable set for the elements defined via the 
above recursions:
\vskip.1in
\noindent
\begin{tabular}{|c|c|c|c|c|}\hline
upper bound on total degree & in $\aa$ & in $\bb$ & in $\AA$ & in $\BB$ \\ \hline\hline
for $\alpha^{(m,n)}$ & $N-1$ & $N-1$ & $N-1$ & $N-1$ \\\hline
for $\beta^{(m,n)}$ & $N$ & $N-1$ & $N$ & $N-1$ \\\hline
for $\gammma_k^{(m,n)}$ & $N$ & $N$ & $N-1$ & $N-1$ \\\hline
for $d$ & $N'$ & $N'$ & $N'$ & $N'-1$ \\\hline
for $e$ & $N''$ & $N''$ & $N''-1$ & $N''$ \\\hline
\end{tabular}

\begin{question}
Is there a good {\it lower} bound on the degree of 
$\alpha, \beta, \gammma_k$ in Proposition~\ref{constructive-Gauss's-lemma}?
\end{question}

\section{Historical remarks}
\label{historical-section}

\subsection{Trivial rings}
The idea from the proof in Section~\ref{less-constructive-section} to show 
that the ring $Q^{(m,n)}$ is trivial is not at all new. 
See Richman \cite{Richman-trivial-ring} for four examples of this idea,
applied in a constructive manner, to prove results in commutative algebra.

\subsection{Lemmas of Artin, Dedekind-Mertens, Gauss-Joyal, 
Kronecker, and McCoy}
Gau\ss's Lemma is closely related to various results by 
Artin, Dedekind-Mertens, Gau\ss-Joyal, Kronecker, McCoy. Proofs of these results avoiding any variants of Zorn's Lemma,
as well as the relations between them and their history, 
are beautifully discussed by 
Coquand, Ducos, Lombardi, and Quitt\'e \cite{CDLQ}, and in
the book by Lombardi and Quitt\'e \cite[\S II Lemme 2.6, \S III.2, \S III.3, Exercice III-6, and Probl\'eme IX-3]{LombardiQuitte}.  

\subsection{Zorn's lemma versus existence of maximal ideals}
As mentioned in Section~1, standard proofs of Gau\ss's Lemma appeal 
not to Zorn's Lemma itself, 
but to the existence, for a proper ideal in a commutative ring, of a
maximal ideal containing it.  
As shown by Hodges \cite{Hodges-KrullimpliesZorn}, already the existence of 
maximal ideals in unique factorization domains implies the axiom of choice.
However, for traditional proofs of Gau\ss's lemma, 
it suffices to know that any proper ideal is contained in a prime ideal.
To explain this,
recall that a {\it Boolean ring} is a ring, necessarily commutative, in which 
every element is idempotent. The fact that non-trivial Boolean rings contain 
prime (and hence maximal) ideals was proven by Stone \cite{Stone}
and is called the Prime Ideal Theorem for Boolean rings. It was announced by 
Scott \cite{Scott} that this implies that any non-trivial commutative ring 
has a prime ideal. Proofs of this result are given in 
Banaschewski \cite{Banaschewski} and Rav \cite{Rav}. 
Perhaps the simplest proof, which the first author learned from O. Gabber, 
is based upon the fact, proven in Olivier \cite{Olivier}, 
that any commutative ring, $R$, has an homomorphism $R \rightarrow T(R)$ 
to a {\it von-Neumann regular} commutative ring, also called an {\it absolutely 
flat ring}, which is universal for homomorphisms to absolutely flat rings. 
The set of idempotents $E(R)$ of $T(R)$ form a 
Boolean ring where the sum of idempotents, $e$ and $e'$ is given by 
$e + e' - 2ee'$. The proof concludes by showing (see Popescu and Vraciu 
\cite{PopescuVraciu}) that the map $\mathfrak p\mapsto\mathfrak p\cap E(R)$ 
gives a  bijection $\Spec(T(R)) \rightarrow \Spec(E(R))$. 
Halpern and L\'evy \cite{HalpernLevy} proved that the
Prime Ideal Theorem is strictly weaker than the axiom of choice.  
Finally note that for non-trivial finitely generated 
$\mathbb Z$-algebras, the existence of maximal ideals does not require 
Zorn's lemma, and is proven in Hodges \cite{Hodges}.

\subsection{The McCoy-Nagata lemma instead of Zorn's}
\label{Nagata-section}

We mention here an alternate proof of 
Proposition~\ref{constructive-Gauss's-lemma}, avoiding Zorn's Lemma,  
which can be deduced from a lemma of McCoy \cite{McCoy}, reproven by Nagata\footnote{We thank the referee for pointing out to us McCoy's paper.  Interestingly,
Nagata states  \cite[p. 213]{Nagata} that "The writer does not know any existing literature which contains (6.13)."  It is unfortunately too late to point 
out this error either to Nagata, who is deceased, 
or to his publisher Interscience, which no longer exists.} 
in \cite[\S 6, pp. 17, 18]{Nagata}.

There Nagata introduces, for a commutative ring $R$,
the set $\SSS$ of all polynomials 
$A(\XX)=\sum_{i=0}^m A_i \XX^i$ in $R[\XX]$ whose coefficients $\{A_i\}_{i=0}^m$
generate the unit ideal $R$.  He constructs the ring
of rational functions $R(\XX)=\SSS^{-1} R[\XX]$, containing $R[\XX]$
as a subring. To do this, he first notes that
\begin{enumerate}
\item[(i)] $\SSS$ is multiplicatively closed, that is,
Proposition~\ref{Gauss's-lemma} above, and
\item[(ii)] $\SSS$ contains no zero-divisors.
\end{enumerate}

\noindent
He offers no proof for assertion (i), but notes
that (ii) is immediate from the following lemma,
which he proves without recourse to Zorn.
\begin{lemma} \cite[Theorem 2]{McCoy}, \cite[(6.13)]{Nagata}
\label{Nagata's-lemma}
For $Q$ a commutative ring, a nonzero element
$A(\XX)$ in $Q[\XX]$ is a zero divisor if and only if there exists 
$q \neq 0$ in $Q$ such that $q A(\XX) = 0$.
\end{lemma}

\noindent
Lemma~\ref{Nagata's-lemma} also leads to a 
proof without Zorn's Lemma of Proposition~\ref{constructive-Gauss's-lemma}
(and hence of Gau\ss's Lemma), as we now explain.

Choose $Q=Q^{(m,n)}$ to be the quotient of $S=\ZZ[\aa,\bb,\AA,\BB]$ 
by the ideal generated by
$a, b, C_0,\ldots,C_{m+n}$, as in Section~\ref{less-constructive-section},
so that Proposition~\ref{constructive-Gauss's-lemma} 
asserts that $Q$ is the zero ring.
Assume for the sake of contradiction that 
$Q$ is {\it not} the zero ring.  Then
$A(\XX)$ and $B(\XX)$ both have nonzero images in $Q[\XX]$, because
their coefficients generate the unit ideal of $Q$.
However, their product $C(\XX)$ has zero image in $Q[\XX]$.
Hence the image of $A(\XX)$ in $Q[\XX]$ is a zero divisor.
By Lemma~\ref{Nagata's-lemma}, there exists $q \neq 0$ in $Q$ 
such that $q A(\XX) = 0$, leading to the contradiction
$q=q \cdot 1 = \sum_{i} a_i q A_i = 0$.

\subsection{Other constructive proofs}
As Richman's \cite[Theorem 4]{Richman-division} is proved
constructively, it immediately gives an algorithmic proof of Gau\ss's Lemma.
It is likely that the same is true for the proof 
of \cite[\S II Lemme 2.6]{LombardiQuitte}, and the McCoy-Nagata proof 
discussed above, at least in the case of ``discrete'' rings\footnote{A ring is said to be {\it discrete} provided that one can constructively decide whether an element of the ring equals zero.}.

\section{Acknowledgements}
The authors thank Kris Fowler, Nick Katz, Henri Lombardi, 
Ezra Miller, Fred Richman, Dennis Stanton, Volkmar Welker,
and the anonymous referee for helpful comments and references.


\begin{thebibliography}{99}

\bibitem{Banaschewski}
B. Banaschewski.
The Power of the Ultrafilter Theorem.
{\it J. London Math. Soc.} {\bf 27} (1983), 193--202

\bibitem{CDLQ}
T. Coquand, L. Ducos, H. Lombardi, C. Quitt\'e.
L'id\'eal des coefficients du produit de deux polyn\^{o}mes. 
{\it Revue des Math\'ematiques de l'Enseignement Sup\'erieur}
{\bf 113} (3), (2003), 25--39. 

\bibitem{HalpernLevy}
J.D.  Halpern and A. Levy.
The Boolean Prime Ideal theorem does not imply the axiom of choice.
Axiomatic Set Theory, {\it Proc. Symp, Pure Math.} {\bf 13}, Vol I, 83 -- 134,
Amer. Math. Society, 1971.

\bibitem{Hodges}
W. Hodges.
Six impossible rings.
{\it J. Algebra} {\bf 31} (1974), 218--244. 

\bibitem{Hodges-KrullimpliesZorn}
W. Hodges.
Krull implies Zorn.
{\it J. London Math. Soc.} {\bf 19} (1979), 285--287.

\bibitem{LombardiQuitte}
H. Lombardi and C. Quitt\'e.
Alg\`ebre commutative: 
M\'ethodes constructives,
Modules projectifs de type fini.
Calvage and Mounet, 2011.


\bibitem{McCoy}
N.H. McCoy,
Remarks on divisors of zero.
{\it Amer. Math. Monthly} {\bf 49} (1942), 286--295.

\bibitem{Nagata}
M. Nagata.
Local rings.
{\it Interscience Tracts in Pure and Applied Mathematics} {\bf 13}. 
Interscience Publishers a division of John Wiley \& Sons,
New York-London 1962

\bibitem{Olivier}
J.-P. Olivier. 
Anneaux absolument universel et epimorphismes \`a but reduit,  Expos\'e VI in
{\it S\'eminaire d`Alg\`ebre Commutative dirig\'e par Pierre Samuel: 1967/1968. 
Les \'epimorphismes d'anneaux}. Secr\'etariat math\'ematique, 
Paris 1968 (Math Review \# 0245561).


\bibitem{PopescuVraciu}
N. Popescu and C. Vraciu.
Sur la structure des anneaux absoluments plats commutatifs.
{\it J. Algebra} {\bf 40} (1976), 364--383. 

\bibitem{Rav}
Y. Rav.
Variants of Rado's selection lemma and their applications.
{\it Math. Nachr.} {\bf 79} (1977), 145 -- 165. 

\bibitem{Richman-trivial-ring}
F. Richman. 
Nontrivial uses of trivial rings.
{\it Proc. Amer. Math. Soc.} {\bf 103} (1988), 1012--1014. 

\bibitem{Richman-division}
F. Richman. 
A division algorithm. 
{\it J. Algebra Appl.} {\bf 4} (2005), no. 4, 441--449. 

\bibitem{Scott}
D. Scott.
Prime ideal theorems for rings, lattices, and Boolean algebras.
{\it Bull. Amer. Math. Soc.} {\bf 60} (1954), 390.

\bibitem{Stone}
M. H. Stone, The theory of representations of boolean algebras.
{\it Trans. Amer. Math. Soc.} {\bf 40} (1936), 37 -- 111.



\end{thebibliography}
\end{document}